\documentclass[10pt,notitlepage,twoside,a4paper]{amsart}
 \usepackage{amsfonts}

\usepackage{amsmath,amssymb,enumerate}

\usepackage{epsfig,fancyhdr,color}

\usepackage{amssymb}
\usepackage{amsmath,amsthm}
\usepackage{latexsym}
\usepackage{amscd}
\usepackage{psfrag}
\usepackage{graphicx}
\usepackage[latin1]{inputenc}
\usepackage[mathcal]{eucal}

\renewcommand{\textsc}{\textcolor{red}}

\newtheorem{theorem}{\rm\bf Theorem}[section]
\newtheorem{proposition}[theorem]{\rm\bf Proposition}
\newtheorem{lemma}[theorem]{\rm\bf Lemma}

\newtheorem*{theorem 1}{\rm\bf Proposition 1}
\newtheorem*{theorem 2}{\rm\bf Proposition 2}

\theoremstyle{definition}

\theoremstyle{remark}
\newtheorem{remark}[theorem]{\rm\bf Remark}

\newtheorem{question}[theorem]{\rm\bf Question}

\def\interieur#1{\mathord{\mathop{\kern 0pt #1}\limits^\circ}}

\title[Lipschitz maps]{Some Lipschitz maps between hyperbolic surfaces with applications to Teichm\"uller theory}

\author{Athanase Papadopoulos}
\address{Athanase Papadopoulos, Max-Planck-Institut f\"ur Mathematik, Vivatsgasse 7, 53111 Bonn, Germany, and : Institut de Recherche Math{\'e}matique Avanc\'ee,
Universit{\'e} de Strasbourg and CNRS,
7 rue Ren\'e Descartes,
 67084 Strasbourg Cedex, France (address for correspondence)} \email{papadopoulos@math.u-strasbg.fr}
\date{\today}

\author{Guillaume Th\'eret}
\address{Guillaume Th\'eret, Max-Planck-Institut f\"ur Mathematik, Vivatsgasse 7, 53111 Bonn, Germany}
\email{guillaume.theret71@orange.fr}

\begin{document}

\begin{abstract}  
In the Teichm\"uller space of a hyperbolic surface of finite type, we construct geodesic lines for Thurston's asymmetric metric having the property that when they are traversed in the reverse direction, they are also geodesic lines (up to reparametrization). 
The lines we construct are special stretch lines in the sense of Thurston. They are directed by complete geodesic laminations that are not chain-recurrent, and they have a nice description in terms of Fenchel-Nielsen coordinates. 
At the basis of the construction are certain maps with controlled Lipschitz constants between right-angled hyperbolic hexagons having three non-consecutive edges of the same size. 
Using these maps, we obtain Lipschitz-minimizing maps between hyperbolic particular pairs of pants and, more generally, between some hyperbolic sufaces of finite type with arbitrary genus and arbitrary number of boundary components. 
The Lipschitz-minimizing maps that we contruct are distinct from Thurston's stretch maps.

\bigskip

\noindent AMS Mathematics Subject Classification:   32G15 ; 30F30 ; 30F60.
\medskip

\noindent Keywords: Teichm\"uller space, surface with boundary, Thurston's asymmetric metric, stretch line, stretch map, geodesic lamination, maximal maximally stretched lamination, Lipschitz metric.
 
\end{abstract}
\maketitle

\section{Introduction}
\label{intro}

In this paper, we prove some results on Thurston's asymmetric metric on Teichm/"uller space. This metric was introduced by Thurston in his paper \ref{Thurston-stretch}.

 we start by constructing Lipschitz homeomorphisms with controlled Lipschitz constant between {\it symmetric right-angled hyperbolic hexagons}, that is, convex right-angled hyperbolic hexagons having three non-adjacent edges of equal length. 
Using these Lipschitz homeomorphisms, we obtain, by doubling the hexagons, Lipschitz homeomorphisms
between {\it symmetric hyperbolic pairs of pants}, that is, hyperbolic pairs of pants which have three geodesic boundary components of equal lengths. 
These Lipschitz homeomorphisms between symmetric pairs of pants are extremal in the sense that their Lipschitz constant is minimal among all Lipschitz constants of homeomorphisms in the same isotopy class.  
But these Lipschitz extremal homeomorphisms between pairs of pants are \emph{not}  stretch maps in the sense of Thurston. 
By varying the Lipschitz constants of the homeomorphisms we construct, we obtain a path in the Teichm\"uller space of the pair of pants which actually coincides with a stretch line in the sense of Thurston, and we exploit the properties of such stretch lines.

We recall that stretch lines are geodesics with respect to Thurston's asymmetric metric, defined  by minimizing the Lipschitz constant between marked hyperbolic surfaces. 
 
By gluing pairs of pants along their boundary components, and by combining the maps we construct between pairs of pants, we obtain stretch lines in the Teichm\"uller space of  hyperbolic surfaces of finite type, of arbitrary genus and of arbitrary number of boundary components, which are also geodesics  (up to reparametrization), for Thurston's asymmetric metric, when they are traversed in the opposite direction. 
These are the first examples we know of such geodesics for this metric.
 
We also recall that by a result of Thurston, given any two points $g$ and $h$ in Teichm\"uller space, there is a unique maximally stretched chain-recurrent geodesic lamination $\mu(g,h)$ from $g$ to $h$ which is maximal (with respect to inclusion), and that if $g$ and $h$ lie in that order on a stretch line directed by a complete chain-recurrent geodesic lamination $\mu$, then $\mu(g,h)=\mu$. 
We obtain the following results that are variations on this theme: We show that if two elements $g$ and $h$ in Teichm\"uller space lie (in that order) on a stretch line we construct, the lamination $\mu(g,h)$ is strictly smaller than the lamination that directs that line, and that there are several (non chain-recurrent) maximal maximally stretched geodesic laminations from $g$ to $h$. 
In other words, the stretch lines we construct are directed by complete geodesic laminations that are {\it not} chain-recurrent, and unlike the chain-recurrent case, these laminations are not uniquely defined.

\section{Thurston's stretch maps between hyperbolic ideal triangles and between pairs of pants}\label{s:stretch}

In this section, we recall the definition of a stretch map between hyperbolic ideal triangles and between pairs of pants. This construction is due to Thurston (see \cite{Thurston-stretch}).

We start with a stretch map from a hyperbolic ideal triangle to itself.  

\begin{figure}[!hbp]
\centering
\psfrag{a}{\small \shortstack{horocycles\\perpendicular\\to the boundary}}
\psfrag{b}{\small \shortstack{horocyclic arc\\of length one}}
\psfrag{c}{\small \shortstack{non-foliated\\region}}
\includegraphics[width=0.5\linewidth]{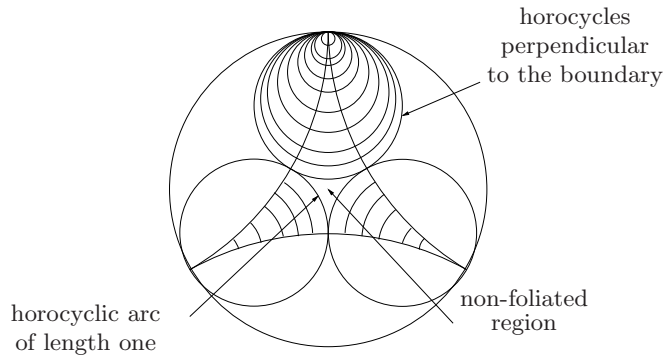}
\caption{\small {The horocyclic foliation of an ideal triangle.}}
\label{horocyclic-foliation}
\end{figure}

Consider a hyperbolic ideal triangle equipped with the partial foliation by horocyclic segments that are perpendicular to the boundary. Up to isometry, there is a unique such object. There is a non-foliated region at the center of the triangle, bounded by three pieces of horocycles (see Figure \ref{horocyclic-foliation}). 
This horocyclic foliation is equipped with a natural transverse measure, which is characterized by the fact that the transverse measure assigned to any arc contained in an edge of the ideal triangle coincides with the Lebesgue measure induced by the hyperbolic metric.  

The non-foliated region of a hyperbolic triangle intersects each edge of the triangle at a point called the {\it center} of that edge.

Let $T$ be the hyperbolic ideal triangle equipped with its horocyclic measured foliation, and consider a real number $k\geq 1$.
The {\it stretch map of magnitude $k$} of $T$ is a homeomorphism $f_k:T\to T$ satisfying the following properties:
\begin{enumerate}
\item The restriction of $f_k$  to the non-foliated region of $T$ is the identity map of that region.
\item On each edge of $T$, $f_k$ sends any point at distance $x$ from the center of that edge to a point at distance $k x$.
\item The map $f_k$ preserves the   horocyclic foliation of $T$; that is, it sends leaves to leaves.
\item On each leaf of the horocyclic foliation, $f_k$ contracts linearly the length of that leaf.
\end{enumerate}

By gluing stretch maps between ideal triangles we construct stretch maps between hyperbolic pairs of pants.

A hyperbolic pair of pants is a sphere with three open disks removed, equipped with a hyperbolic metric in which the three boundary components are closed geodesics (the lift of such a curve to the hyperbolic universal cover seen as a subset of the hyperbolic plane $\mathbb{H}^2$  is a geodesic in $\mathbb{H}^2$).

Let $P$ be a hyperbolic pair of pants. 
We choose a complete geodesic lamination $\lambda$ in $P$. 
Such a complete geodesic lamination necessarily consists of  three disjoint bi-infinite geodesics that spiral around the boundary components of $P$, decomposing that surface into two hyperbolic ideal triangles. 
The horocyclic measured  foliations of the two ideal triangles fit together smoothly since they are both perpendicular to the edges of the ideal triangles, and therefore they form a Lipschitz line field on the surface. 
For each $k\geq 1$, consider a stretch map of magnitude $k$ defined on each of the ideal triangles composing $P$. 
We obtain a new hyperbolic pair of pants $P_k$  by gluing the ideal triangles together along their boundaries according to identifications that are compatible with the stretch maps. 
This defines a homeomorphism from $P$ to another hyperbolic pair of pants $P_{k}$, which is called a {\it stretch map} 
(of magnitude $k$) from $P$ to $P_{k}$.

The above construction can be repeated on several copies of hyperbolic pairs of pants.
By gluing together these pairs of pants according to the identifications given by the stretch maps, 
we obtain a stretch map of magnitude $k$ from a hyperbolic surface $S$ to another $S_{k}$.
Note that the complete geodesic laminations giving the decompositions into ideal triangles of the pairs of pants in $S$ 
give, together with the pants decomposition of $S$, a complete geodesic lamination on the surface $S$.

\begin{remark}
The reader should be aware that stretch maps are actually defined in a much wider generality than the one presented here.
The underlying complete geodesic lamination giving the decomposition of the surface into ideal triangles can be chosen arbitrarily among the complete geodesic laminations and it is not necessarily the completion of a geodesic pants decomposition as above. However, in this paper, we shall only need the special case of stretch maps described above.
\end{remark}

\section{Extremal Lipschitz maps between symmetric right-angled hexagons}\label{s:extremal}

Given two metric spaces $(X,d_X)$ and $(Y,d_Y)$ and a map $f:X\to Y$ between them,  the {\it Lipschitz constant} $\mathrm{Lip}(f)$ of $f$ is defined  as

 \begin{displaymath}\label{Lip}
\hbox{Lip}(f)=\sup_{x\neq y\in X}\frac{d_{Y}\big{(}f(x),f(y)\big{)}}{d_{X}\big{(}x,y\big{)}}\in   \mathbb{R}\cup\{\infty\}.
\end{displaymath}

We shall say that the map $f$ is {\it Lipschitz} if its Lipschitz constant is finite.

The stretch maps $f_k$ between hyperbolic ideal triangles that we considered in the last section are examples of Lipschitz homeomorphisms, with Lipschitz constant equal to $k$.  
Note that the fact that this Lipschitz constant is at least $k$ can be seen from the action of these maps on the boundary of the ideal triangles. 
The fact that the Lipschitz constant is exactly $k$ is implicit in Thurston's paper \cite {Thurston-stretch}. 
It also follows from the computations below (see Remark \ref{rem:Lip}). 
By using these maps as building blocks, we recalled in \S\ref{s:stretch} how one obtains Lipschitz homeomorphisms of hyperbolic pairs of pants and, more generally,
of hyperbolic surfaces. 
These stretch maps have Lipschitz constants $k$.

In this section, we shall define Lipschitz maps between some particular hyperbolic right-angled hexagons, which will also have controlled Lipschitz constants, and which can be used to define Lipschitz homeomorphisms between special hyperbolic pairs of pants, by gluing hyperbolic right-angled hexagons and taking the union of Lipschitz maps between them.
By gluing together these special pairs of pants in an appropriate manner, this will eventually yield homeomorphisms between special hyperbolic surfaces of arbitrary finite type, with controlled Lipschitz constants.

 A {\it symmetric right-angled hexagon} is a geodesic hexagon $H$ in the hyperbolic plane $\mathbb{H}^2$ with three pairwise non-consecutive edges having the same length. 
(Note that this implies that the remaining three edges also have the same length.)

We consider a symmetric right-angled hexagon $H$, and we choose three pairwise non-consective edges of $H$, which we call the {\it long} edges. 
We denote their common length  by $2L$. 
The other three non-consecutive edges are called {\it short}, and we denote their common length by $2l$.
An easy computation using well-known formulae for right-angled hexagons gives
\begin{equation}\label{lengths}
2\sinh(l)\sinh(L)=1.
\end{equation}

For each real number $k\geq 1$, we let $H_{k}$ be the symmetric right-angled hexagon
obtained by multiplying the lengths of the long edges of $H$ by the factor $k$. 
We note that this property determines the isometry type of $H_k$ in a unique way.
We call the edges of $H_k$ that are the images of the long edges of $H$ by this dilatation map
 the {\it long} edges of $H_k$ and we denote their common length by $2L_k$.
We let $2l_{k}$ denote the length of the other edges of $H_{k}$, which we call the {\it short} ones.

In this section, all the maps between symmetric right-angled hexagons that we shall consider will be homeomorphisms sending the long (respectively short) edges  to the long (respectively short) edges, and in general we shall not repeat this condition.

\begin{figure}[!hbp]
\psfrag{A}{$A$}
\psfrag{B}{$B$}
\centering
\includegraphics[width=0.4\linewidth]{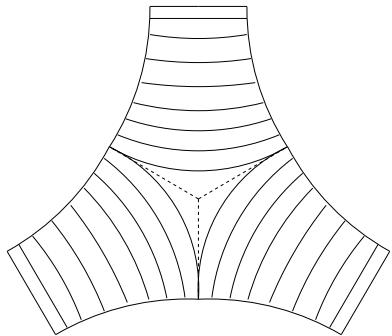}
\caption{\small{The foliation by curves  equidistant to the short edges of a symmetric right-angled hexagon. The central region is not foliated, and it is bounded by three hypercycles which meet each other tangentially.}}
\label{foliated_hexagon}
\end{figure}

The three lengths of any three non-consecutive edges of $H$ (respectively of $H_{k}$) satisfy the triangle inequality.
Therefore, we can equip $H$ (respectively $H_{k}$) with a partial measured foliations $F$ (respectively $F_{k}$) whose leaves are loci of equidistant points from the  short edges. 
In the hyperbolic plane, equidistant points from geodesics are classicaly called {\it hypercycles}, and we shall use this terminology. 
The foliations of $H$ (respectively $H_k$) by hypercycles are shown in Figure \ref{foliated_hexagon}, and such foliations have already been considered by Thurston in his compactification theory of Teichm\"uller space (see \cite[expos\'e 6]{FLP}). 
There is a non-foliated region of $F$ (respectively $F_k$) at the center of $H$ (respectively $H_k$).

The intersection number of $F$ (respectively, $F_{k}$) with an edge of $H$ (respectively, $H_{k}$) is
 either $2L$ or $0$ (respectively, $2kL$ or $0$) depending on whether the edge is long or short.
 
We also equip $H$ (respectively $H_{k}$) with the partial foliation $G$ (respectively $G_{k}$) whose leaves are geodesic arcs 
perpendicular to the leaves of $F$ (respectively $F_{k}$).

In Theorem \ref{theorem:Lipschitz}, we shall construct a map, $h_k: H\to H_{k}$ which (leafwise) sends $F$ to $F_k$, and $G$ to $G_k$ and whose Lipschitz constant is $k$.
Such a map is Lipschitz-extremal in its homotopy class relative to the boundary, since the Lipschitz constant of any map $f:H\to H_k$ which sends long (respectively short) edges of $H$ to long (respectively short) edges of $H_k$ is bounded below by $k$. 
The Lipschitz-extremal maps we shall construct are ``canonical" in the sense that they preserve a pair of hypercyclic/geodesic foliations, and they are reminiscent of Thurston's stretch maps between ideal triangles. 
In some precise sense that we specify below, Thurston's stretch maps between ideal triangles are limits of the Lipschitz-extremal maps between symmetric hexagons.

Before defining the map $h_k$, we make a geometrical remark.
Consider the family of all symmetric right-angled hexagons $H_{k}$ as $k$ varies from $1$ to infinity.
Each of these hexagons has a center which is the center of the rotation that permutes  each triple of non-consecutive edges. 
For each such hexagon, consider the three geodesic rays emanating from its center and meeting the short edges perpendicularly.
Place all the hexagons $H_{k}$ in the hyperbolic plane so that all their centers coincide and such that all the above geodesic rays coincide as well.
Now for each such hexagon $H_{k}$, consider the associated extended hexagon $\widehat{H}_{k}$ defined as the region of infinite area enclosed by the three geodesics in $\mathbb{H}^2$
extending the long edges of $H_{k}$.
It follows from Equation (\ref{lengths}) that as $L_k$ decreases, $l_k$ increases, and conversely. From this, we deduce that for any $1\leq k\leq k'$, we have $\widehat{H}_{k'}\subset \widehat{H}_{k}$.

\begin{figure}[!hbp]
\psfrag{a}{$a$}
\psfrag{b}{$b$}
\psfrag{m}{$a'$}
\psfrag{n}{$b'$}
\psfrag{p}{$p$}
\psfrag{q}{$q$}
\psfrag{O}{$O$}
\centering
\includegraphics[width=0.6\linewidth]{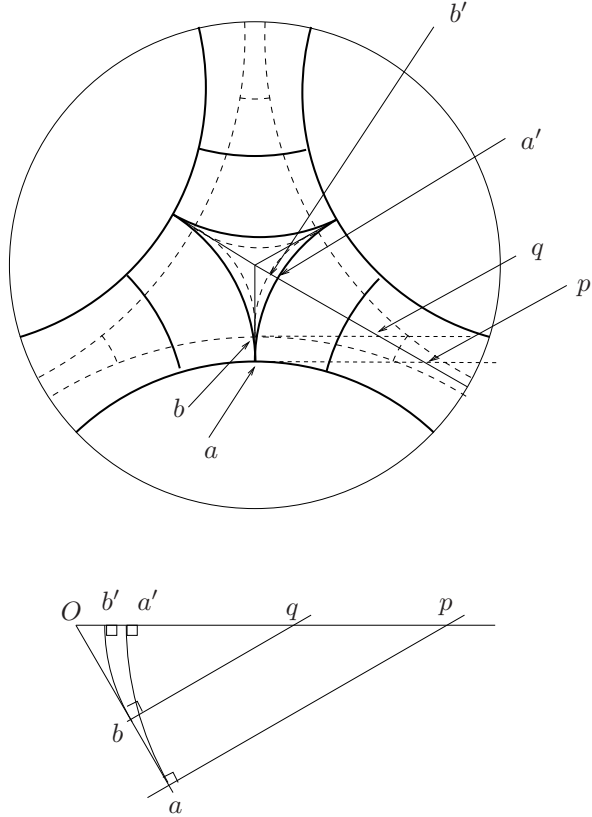}
\caption{\small{The upper figure represents, 
in bold lines, a symmetric right-angled hexagon $H_k$, and in dashed lines, a symmetric right-angled hexagon $H_{k'}$ with $k'>k$, together with their extensions $\hat{H}_k$ and $\hat{H}_k'$. 
The fact that the non-foliated region of the symmetric hexagon $H_k'$ is included in the non-foliated region of the symmetric hexagon $H_k$,  for $k'>k$, as it is represented in the upper figure, can be deduced from the Euclidean construction in the lower figure, in which
 the arcs $aa'$ and $bb'$ are on the boundaries of the non-foliated regions of $H_k$ and $H_{k'}$ respectively.}}
\label{hex1}
\end{figure}

We also note that as $k$ tends to infinity, the extended hexagon $\widehat{H}_{k}$ as well as the hexagon $H_{k}$ itself converge, in the Hausdorff topology associated to the Euclidean metric (using as in Figure \ref{hex1} the disk-model of the hyperbolic plane) to an ideal triangle. Likewise, as $k\to\infty$, 
the measured foliation $F_{k}$ converges to the horocyclic foliation of the ideal triangle (represented in Figure \ref{horocyclic-foliation})
 and the non-foliated region of $F_{k}$ converges to
the non-foliated region of that horocyclic foliation.

The following two lemmas will be used in the proof of Theorem \ref{theorem:Lipschitz} below.

\begin{lemma}\label{lemma:inclusion} 
For $k'>k\geq1$, the non-foliated region of $F_{k'}$ is strictly contained in the non-foliated region of $F_{k}$.
\end{lemma}

\begin{proof} 
We work in the disk model of the hyperbolic plane. 
The statement will follow from the construction of the symmetric hexagons, represented in Figure \ref{hex1}. 
In the upper part of that figure, the hexagon $H_k$ (also with its edges extended) is drawn in bold lines, and the hexagon $H_{k'}$  (with its edges extended) is drawn in dashed lines. 
We have chosen the hexagons to be symmetric with respect to the Euclidean center $O$ of the unit disk. 
In the upper figure, the point $p$ (respectively $q$) is the Euclidean center of the hypercycle that is on the boundary of non-foliated region of $H_k$ (respectively $H_{k'}$). 
The point $a$ (respectively $b$) is a vertex of the non-foliated region of $F_{k}$ (respectively $F_{k'}$). 
A more detailed view of a region drawn in the the upper part of Figure \ref{hex1}  is represented in the lower part. 
The point $a'$ (respectively $b'$) is the center of a boundary hypercycle of the non-foliated region of $F_{k}$ (respectively $F_{k'}$). 
The Euclidean triangles $Opa$ and $Oqb$ are homothetic by a Euclidean homothety of center $O$ and factor $<1$. 
This homothety sends the Euclidean circle arc $aa'$ to the Euclidean circle arc $bb'$. 
Thus, there exists a Euclidean homothety of center $O$ that sends the non-foliated region of  $H_{k'}$ strictly into the non-foliated region of $H_k$, which proves the lemma. 
\end{proof}

\begin{lemma}\label{lemma:computation} 
In the upper half-plane model of the hyperbolic plane, consider the geodesic represented by the imaginary axis $i\mathbb{R}^+= \{ir, r>0\}$, and a hypercycle making an angle $\frac{\pi}{2}-\theta_1$ with this geodesic, with $0< \theta_1<\pi/2$. 
Let $\ell$ be the length of a geodesic arc $\alpha$ joining perpendiculary the vertical geodesic and the hypercycle. 
Then, we have
\[\cos \theta_1 = \tanh \ell.\]
\end{lemma}
\begin{proof}

\begin{figure}[!hbp]
\psfrag{i}{$i$}
\psfrag{l}{$\ell$}
\psfrag{t}{$\theta_1$}
\psfrag{a}{$\alpha$}
\centering
\includegraphics[width=0.35\linewidth]{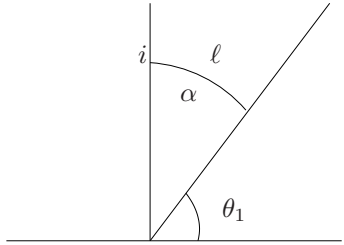}
\caption{\small{$\ell$ is the length of a segment $\alpha$ joining perpendicularly the vertical geodesic and the hypercycle making an angle $\theta_1$ with the horizontal. 
We have $\cos\theta_1=\tanh \ell$.}}
\label{angle-length}
\end{figure}

 We refer to Figure \ref{angle-length}. We parametrize the geodesic arc $\alpha$ by the map
 \[\alpha:[\theta_1,\pi/2]\to \mathbb{H}^2\]
\[\theta\mapsto (\cos\theta, \sin\theta).\]
 Using the formula for the infinitesimal length element in the upper half-plane model, we can write 
 \[\ell = \int_{\theta_{1}}^{\pi/2} \frac {\Vert \alpha'(\theta)\Vert}{\mathrm{Im}(\alpha(\theta))}d\theta = \int_{\theta_{1}}^{\pi/2}\frac{d\theta}{\sin \theta}.\]
 Computing the integral, we find \[e^{-\ell}=\tan(\theta_1/2)\]
 and after transformation we obtain
 \[\cos \theta_1 = \tanh \ell.\]
\end{proof}

 We now construct the map $h_k:H\to H_k$.

From the inclusion of the non-foliatied region of $H_k$ into the non-foliated region of 
$H$ for all $k\geq 1$ (Lemma \ref{lemma:inclusion}), it will follow that the map $h_k$ we shall construct can be chosen to be contracting from the non-foliated region of $H$ to the 
non-foliated region of $H_{k}$. 

To define the map $h_k$, it suffices to do it in a component of the foliated region of $H$.
Consider such a component.
It is isometric to the region $C$ in the upper half-plane model of the hyperbolic plane defined in polar coordinates by
$$
C=\{z=Re^{i\theta}\ :\ 1\leq R\leq e^{2l},\ \theta_{1}\leq\theta\leq\pi/2\},
$$
where $\theta_{1}$ is chosen so that the geodesic parameterized by $\theta\mapsto Re^{i\theta}$, $\theta_{1}\leq\theta\leq\pi/2$,
has length $L$.

From Lemma \ref{lemma:computation}, we have
$$
\cos \theta_{1}=\tanh L.
$$
Likewise, the image by $h_k$ of the component $C$ of the complement in $H$ of the non-foliated region is isometric to 
 the region $C_{k}$ in the upper half-plane model of $\mathbb{H}^2$ given by
$$
C_{k}=\{z=Re^{i\theta}\ :\ 1\leq R\leq e^{2l_{k}},\ \theta_{k}\leq\theta\leq\pi/2\},
$$
where 
$$
\cos(\theta_{k})=\tanh(kL).
$$
In these descriptions, the foliations $F$ and $F_{k}$, are given by the hypercycles defined by $\theta=\mathrm{cst}$, 
while  the foliations $G$ and $G_{k}$, are given by the geodesics defined by $R=\mathrm{cst}$.
The short sides of $C$ and $C_{k}$ correspond to $\theta=\pi/2$.
Our map $h_k$ maps a point $A\in C$ which is at distance $d$ from the short side of $C$ to a point which is at distance $kd$ from the short side of $C_{k}$.
If the point $A$ lies on the leaf of $G$ which cuts the short side of $C$ at distance $h$, 
then the image of $A$ by $h_k$ belongs to the leaf that cuts the short side of $C_{k}$ at distance $hl_{k}/l$.

We need to have an explicit formula for $h_k$ in order to compute the norm of its derivative.

Let $A$ be a point in $C$ given in polar coordinates by $(R,\theta)$.
Denote the coordinates of the point $h_k(A)\in C_{k}$ by $(R',\theta')$.
We also describe the points $A$ and $h_k(A)$ by their distances from the short sides,
namely $d$ and $kd$, and by their distances from the lowest geodesic boundary   
of $C$ and $C_{k}$, as above.
 
Let us first compute $R'$.
The logarithm of $R$ and of $R'$ are the distances of the points $A$ and $h_k(A)$
from  the lowest geodesic boundary of $C$ and $C_{k}$, respectively.
By what has been previously said,
we have
$$
\log R'=\frac{l_{k}}{l}\log R.
$$
Therefore,
$$
R'=R^{l_{k}/l}.
$$

Let us now compute $\theta'$.
The same computation as for the formula giving $\theta_{1}$ establishes
$$
\sin \theta=\frac{1}{\cosh d },\ \textrm{or}\  \cos \theta =\tanh d.
$$
Therefore,
$$
d=\textrm{argcosh}\Big{(}\frac{1}{\sin \theta}\Big{)}.
$$
Now,
$$
\theta'=\textrm{arccos}(\tanh(kd)).
$$

Thus we get the following formula for $h_k$, viewed as a map from $C$ to $C_{k}$,
$$
h_k(R,\theta)=\left(R^{l_{k}/l},\textrm{arccos}(\tanh\big{(}k\,\textrm{argcosh}\Big{(}\frac{1}{\sin \theta}\Big{)}\big{)})\right).
$$


Now that the homeomorphism $h_k$ is defined, we proceed to show that its Lipschitz constant equals $k$. 
For this, we compute the norm of its derivative.

We easily have
$$
\frac{\partial R'}{\partial R}=\frac{l_{k}}{l}R^{(l_{k}/l)-1},\ \frac{\partial R'}{\partial \theta}=0,\ \frac{\partial \theta'}{\partial R}=0.
$$

Since $\textrm{arccos}'(x)=\displaystyle -\frac{1}{\sqrt{1-x^2}}$, we get
\begin{eqnarray*}
\frac{\partial\theta'}{\partial\theta}&=&-\frac{1}{\sqrt{1-\tanh^{2}(k\,\textrm{argcosh}\Big{(}\frac{1}{\sin \theta }\Big{)})}}
\frac{\partial}{\partial\theta}\Big{(}\tanh\big{(}k\,\textrm{argcosh}\Big{(}\frac{1}{\sin \theta }\Big{)}\big{)}\Big{)}\\
&=&-\cosh(k\,\textrm{argcosh}\Big{(}\frac{1}{\sin \theta}\Big{)})
\frac{\partial}{\partial\theta}\Big{(}\tanh\big{(}k\,\textrm{argcosh}\Big{(}\frac{1}{\sin \theta}\Big{)}\big{)}\Big{)}.
\end{eqnarray*}

Now, since $\tanh'(x)=\displaystyle \frac{1}{\cosh^{2}(x)}$, we have
$$
\frac{\partial}{\partial\theta}\Big{(}\tanh\big{(}k\,\textrm{argcosh}\Big{(}\frac{1}{\sin \theta}\Big{)}\big{)}\Big{)}=
\frac{k}{\cosh^{2}(k\,\textrm{argcosh}\Big{(}\frac{1}{\sin \theta}\Big{)})}
\frac{\partial}{\partial\theta}\textrm{argcosh}\Big{(}\frac{1}{\sin \theta}\Big{)}.
$$
Hence, since $\textrm{argcosh}'(x)=\frac{1}{\sqrt{x^2-1}}$,
\begin{eqnarray*}
\frac{\partial\theta'}{\partial\theta}&=&
\frac{-k}{\cosh(k\,\textrm{argcosh}\Big{(}\frac{1}{\sin \theta}\Big{)})}
\frac{\partial}{\partial\theta}\textrm{argcosh}\Big{(}\frac{1}{\sin \theta}\Big{)}\\
&=&\frac{-k}{\cosh(k\,\textrm{argcosh}\Big{(}\frac{1}{\sin \theta}\Big{)})}
\frac{1}{\sqrt{\frac{1}{\sin^{2}\theta}-1}}
\frac{\partial}{\partial\theta}\frac{1}{\sin \theta}\\
&=&\frac{k\sin \theta}{\cos \theta \cosh(k\,\textrm{argcosh}\Big{(}\frac{1}{\sin \theta}\Big{)})}\frac{\cos \theta }{\sin^2 \theta }.
\end{eqnarray*}
Finally, we have
$$
\frac{\partial\theta'}{\partial\theta}=
\frac{k}{\sin \theta }\Big{[}\cosh(k\,\textrm{argcosh}\Big{(}\frac{1}{\sin \theta }\Big{)})\Big{]}^{-1}.
$$
The last partial derivative can also be written as
$$
\frac{\partial\theta'}{\partial\theta}=k\frac{\cosh d}{\cosh(kd)}.
$$

We now proceed to compute the norm of the differential $dh_k$.
Recall that the square of the norm of a vector $(dx,dy)$ in the tangent plane $T_{z}(\mathbb{H}^{2})$ of the upper half-plane model of the hyperbolic plane is given by
$$
\frac{dx^{2}+dy^{2}}{y^{2}},
$$
where $z=x+iy$.
In polar coordinates, this is written as 
$$
\frac{dR^{2}+R^{2}d\theta^{2}}{R^{2}\sin^{2}\theta}.
$$
Let $V=(V_{R},V_{\theta})$ be a non-zero tangent vector at the point $(R,\theta)$.
We compute the norm of the differential $dh_k$ at the point $(R,\theta)$. We have
\begin{eqnarray*}
||(dh_k)_{(R,\theta)}\cdot V||^{2}&=&||(\frac{\partial h_{k}}{\partial R}dR+\frac{\partial h_{k}}{\partial \theta}d\theta)\cdot V||^{2}\\
&=&\frac{1}{R^{2}\sin^{2}\theta}
\Big{(}\big{(}\frac{\partial R'}{\partial R}V_{R}+\frac{\partial R'}{\partial \theta}V_{\theta}\big{)}^{2}+
R^{2}\big{(}\frac{\partial \theta'}{\partial R}V_{R}+\frac{\partial \theta'}{\partial \theta}V_{\theta}\big{)}^{2}\Big{)}\\
&=&\frac{1}{R^{2}\sin^{2}\theta}
\Big{(}\big{(}\frac{\partial R'}{\partial R}V_{R}\big{)}^{2}+
R^{2}\big{(}\frac{\partial \theta'}{\partial \theta}V_{\theta}\big{)}^{2}\Big{)}.
\end{eqnarray*}
Note that
$$
||V||^{2}=\frac{1}{R^{2}\sin^{2}\theta}(V_{R}^{2}+R^{2}V_{\theta}^{2}).
$$
Therefore, since $||(dh_k)_{(R,\theta)}||=\sup_{V\neq0}\frac{||(dh_{k})_{(R,\theta)}\cdot V||}{||V||}$, we get
\begin{eqnarray*}
||(dh_k)_{(R,\theta)}||^{2}&=&\sup_{V\neq0}\Big{(}
\frac{\big{(}\frac{\partial R'}{\partial R}V_{R}\big{)}^{2}+
R^{2}\big{(}\frac{\partial \theta'}{\partial \theta}V_{\theta}\big{)}^{2}}{V_{R}^{2}+R^{2}V_{\theta}^{2}}
\Big{)}\\
&=&\sup_{V\neq0}\Big{(}
\frac{\big{(}\frac{\partial R'}{\partial R}V_{R}\big{)}^{2}+
\big{(}\frac{\partial \theta'}{\partial \theta}R\,V_{\theta}\big{)}^{2}}{V_{R}^{2}+(R\,V_{\theta})^{2}}
\Big{)}\\
&=&\sup_{V_{R}^{2}+(R\,V_{\theta})^{2}=1}\Big{(}
\big{(}\frac{\partial R'}{\partial R}V_{R}\big{)}^{2}+
\big{(}\frac{\partial \theta'}{\partial \theta}R\,V_{\theta}\big{)}^{2}
\Big{)}\\
&=&\max\Big{\{}\Big{(}\frac{\partial R'}{\partial R}\Big{)}^{2},\Big{(}\frac{\partial \theta'}{\partial \theta}\Big{)}^{2}\Big{\}}.
\end{eqnarray*}

We have 
$$
1\leq R\leq e^{2l}.
$$
Since $l_{k}/l\leq 1$, we get
$$
1\geq R^{l_{k}/l-1}\geq e^{2(l_{k}-l)}>0,
$$
that is,
$$
0\leq\frac{\partial R'}{\partial R}\leq 1.
$$
Now, since
$$
\frac{\partial\theta'}{\partial\theta}=k\frac{\cosh(d)}{\cosh(kd)},
$$
we get, for all $(R,\theta)$,
$$
0\leq\frac{\partial\theta'}{\partial\theta}\leq k
$$ 
and the equality $\frac{\partial\theta'}{\partial\theta}=k$ is realized at the points $d=0$, 
that is, on the short side of $C$.
Therefore, we obtain
$$
\sup_{(R,\theta)\in C}||(dh_k)_{(R,\theta)}||=k.
$$

The supremum of the norm of $dh_k$ bounds from above the Lipschitz constant of $h_k$:
If $x,y$ are two points of $C$ and if $\gamma$ is the geodesic path from $x$ to $y$,
we get
\[d(h_k(x),h_k(y))\leq l(h_k(\gamma))=\int_{0}^{d(x,y)}||(dh_k)_{\gamma(t)}\cdot\gamma'(t)||dt\leq\sup_{z}||(dh_{k})_{z}||d(x,y).\]
Therefore, if $L(h_k)$ denotes the Lipschitz constant of $h_k$, we get from what precedes,
$$
L(h_k)\leq k.
$$
Since the long edges are dilated by the factor $k$, we have $L(h_k)\geq k$.
Finally, 
$$
L(h_k)=k.
$$
Putting all pieces together, the map we constructed from $H$ to $H_{k}$ has Lipschitz constant $k$.

We summarize the preceding construction in the following: 

\begin{theorem}\label{theorem:Lipschitz}  
The map $h_k:H\to H_k$ is $k$-Lipschitz.
Furthermore for any $k'<k$, there is no $k'$-Lipschitz map from $H$ to $H_k$. 
\end{theorem}

\begin{proof}
The first part follows from the construction.
Since, by definition, a map $h_k:H\to H_k$ sends the long edges of $H$ to the long edges of $H_k$, we immediately get $\mathrm{Lip}(h_k)\geq k$.  
This proves the second part of the theorem.
\end{proof}

\begin{remark} \label{rem:Lip} 
We already observed that, 
reasoning in the disk model of the hyperbolic plane and using the notion of Hausdorff convergence on bounded closed subsets of that disk with respect to the underlying Euclidean metric, we can make a sequence of symmetric right-angled hexagons converge to an hyperbolic ideal triangle, in such a way that the following three properties hold:

(1) The partial measured foliation of the hexagons by hypercycles converges to the partial measured foliation of the hyperbolic ideal triangle by horocycles.

(2) The partial foliation of the hexagons by geodesics perpendicular to the foliation by hypercycles converges to the partial foliation of the ideal triangle by geodesics perpendicular to the horocycles.

(3) The non-foliated regions of the hexagons converge to the non-foliated region of the ideal triangle.

Furthermore, for all $k\geq 1$, we can make the convergence of hexagons to the ideal triangle in such a way that $k$-Lipschitz maps $f_k:H\to H_k$ converge uniformly on compact sets to the stretch maps $f_k:T\to T$ between hyperbolic ideal triangles. 
This shows in particular that the stretch maps $f_k$ have Lipschitz constant $k$.
\end{remark}

We note that Lipschitz maps between pairs of pants are also considered by Otal in his paper \cite{Otal},  in relation with the Weil-Petersson metric of Teichm\"uller space.

\section{Asymmetric metrics on Teichm\"uller spaces of surfaces with or without boundary}

In this section, $S$ is a surface of finite type $(g,b)$, which may have empty or nonempty boundary ($g$ denotes the genus of $S$ and $b$ the number of boundary components). We assume that the Euler characteristic of $S$ is negative.
The hyperbolic structures we construct on $S$ are such that all the boundary components are closed smooth geodesics. 
We denote by $\mathcal{T}(S)$ or by $\mathcal{T}_{g,b}$ the Teichm\"uller space of $S$, that is, the space of homotopy classes of hyperbolic metrics on that surface. 

Given two hyperbolic structures $X$ and $Y$ on $S$, we define
\begin{equation}\label{def:L}
L(X,Y)=\log \inf_f \mathrm{Lip}(f)
\end{equation}
where the infimum is taken over the set of  Lipschitz homeomorphisms $f:X\to Y$ that are homotopic to the identity.

\begin{lemma}[Thurston] \label{lem:separates}
For any two hyperbolic metrics $X$ and $Y$ on $S$, if
 $L(X,Y)\leq 0$, then $X$ and $Y$ are isometric by a homeomorphism that is homotopic to the identity.
\end{lemma}

\begin{proof}
We follow Thurston's proof of the corresponding result in the case of surfaces without boundary, cf. \cite[Proposition 2.1]{Thurston-stretch}.
Since $L(X,Y)\leq 0$, there exists a sequence of homeomorphisms $f_n:X\to Y$, $n =0,1,\ldots$, with Lipschitz constants  $\mathrm{Lip}(f_n)$ converging to a real number $L\leq 1$. The sequence $(f_n)$ is uniformly equicontinuous, therefore up to taking a subsequence, we can assume that $(f_n)$ converges uniformly to a map $f:X\to Y$. 
We have  $\mathrm{Lip}(f)=L\leq 0$. 
We now prove that $f$ is surjective. 
Take a point $y$ in $Y$, and for all $n\geq 0$, let $x_n=f_n^{-1}(y)$. 
Up to taking a subsequence of $(f_n)$, we can assume, by compactness, that $x_n\to x\in X$. 
We show that $f(x)=y$. 
Let us fix some $\epsilon >0$. 
We have 
\[\vert f(x)-y\vert = \vert f(x)-f_n(x_n)\vert \leq \vert f(x)-f_n(x)\vert + \vert f_n(x) - f_n(x_n)\vert.\]

Since $f_n\to f$ uniformly, there exists $N\geq0$ such that for all $n\geq N$, we have $\vert f(x)-f_n(x)\vert \leq \epsilon/2$.
Since the family $(f_n)$ is equicontinuous, there exists $\delta >0$ such that for $x_1$ and $x_2$ satisfying $\vert x_1-x_2\vert <\delta$, we have $\vert f_m(x_1)-f_m(x_2)\vert \leq \epsilon/2$ for all $m\geq 0$.

Since $x_n\to x$, there exists $N'$ such that for all $n\geq N'$, we have $\vert x-x_n\vert <\delta$.

For $n\geq \max\{N,N'\}$, we have, for all $m$, $\vert f_m(x)-f_m(x_n)\vert \leq \epsilon/2$. 
In particular, for $m=n$, $\vert f_n(x)-f_n(x_n)\vert \leq \epsilon/2$.
This shows that for every $\epsilon >0$, we have $\vert f(x)-y\vert \leq \epsilon$. 
Thus, $f(x)=y$. 
This shows that $f$ is surjective.

We cover $S$ by a set of geometric disks with disjoint interior whose total area is equal to the area of $X$. 
The metrics $X$ and $Y$ have the same area. 
Since $\mathrm{Lip}(f)\leq 1$ and since $f$ is surjective, the image by $f$ of a disk of radius $R$ is a disk of radius $R$. 
Furthermore, $f$ sends the boundary of any such disk to the boundary of the image disk.
We deduce that any geometric disk is sent by $f$ isometrically to a geometric disk of the same radius.
Furthermore, it is easy to see that the center of such a disk is sent to the center of the image disk.

From this, we deduce that $f$ is locally distance-preserving. 
This implies that $f$ is an isometry.
\end{proof}

We call an {\it asymmetric metric} on a set $X$ a function that satisfies the axioms of a metric except the symmetry axiom, and that does not satisfy this axiom.

\begin{proposition}
The function $L$ defined in (\ref{def:L}) is an asymmetric metric on the Teichm\"uller space $\mathcal{T}(S)$.
\end{proposition}

\begin{proof}
By Lemma \ref{lem:separates}, $L$ is nonnegative and separates points. The
triangle inequality is obviously satisfied. The fact that the metric does not
satisfy the symmetry axiom can be seen using an example analogous to the one showing the
corresponding result for surfaces without boundary, given by Thurston in \cite{Thurston-stretch}.
\end{proof}

We let $\mathcal{S}$ be the set of isotopy classes of simple closed curves on $\mathcal{S}$ which are not homotopic to a point (the boundary components of $S$ are included).

The asymmetric metric $L$ is an analogue, for surfaces with boundary, of the asymmetric metric defined by Thurston in \cite{Thurston-stretch} for surfaces without boundary. 
In the same paper, Thurston defined the following function on the Teichm\"uller space $\mathcal{T}(S)$ of a surface $S$ without boundary:
\begin{equation}\label{def:K}
K(x,y)=\log\sup_{C\in\mathcal{S}}\frac{l_y(C)}{l_x(C)}.
\end{equation}

Thurston proved that we obtain the same function $K$ if instead of taking the infimum over the elements of $\mathcal{S}$ in (\ref{def:K}) we take the infimum over all (not necessarily simple) closed curves (see \cite{Thurston-stretch}, Proposition 3.5).

In the case where the surface $S$ has nonempty boundary, Formula (\ref{def:K}) does not define an asymmetric metric on the Teichm\"uller space of $S$. 
This can easily be seen in the case where the surface is a pair of pants $P$. 
Denoting by $C_1,C_2,C_3$ the three boundary components of the pair of pants, the function $K$ defined on $\mathcal{T}(P)\times \mathcal{T}(P)$ takes the form
\[K(x,y)=\log\sup_{i=1,2,3}\frac{l_y(C_i)}{l_x(C_i)}.\]
This function $K$ on $\mathcal{T}(P)$ satisfies the triangle inequality, but it is not an asymmetric metric, since it can take negative values. Furthermore, it does not separate points; that is, there exist distinct $x$ and $y$ in $\mathcal{T}(P)$ with $K(x,y)=0$ (take $x$ and $y$ satisfying $l_x(C_1)=l_y(C_1)$, and $l_x(C_i)>l_y(C_i)$ for $i=2,3$).

In fact, for any surface $S$ with nonempty boundary, there exist hyperbolic metrics $X$ and $Y$ such that $K(x,y)<0$ (see \cite{PT3}).

We have $K\leq L$. Indeed, for any $k$-Lipschitz homeomorphism from a hyperbolic metric $x$ on $S$ to a hyperbolic metric $y$ on $S$, we easily see that we have, for every simple closed curve $\gamma$ on $S$, $l_{y}(f(\gamma))\leq kl_{x}(\gamma)$, which implies $K(x,y)\leq L(x,y)$. 

There is a modification of the function $K$ defined  in Formula (\ref{def:K}) which is adapted to the case of surfaces with or without boundary, which we studied in \cite{LPST} and which we now recall. 
The definition involves considering essential arcs in $S$ together with essential simple closed curves. 
We call an {\it essential arc} in $S$ an embedding of a closed interval, the arc having its endpoints on the boundary of $S$ and its interior in the interior of $S$, and such that this arc is not homotopic relative endpoints to an arc contained  in $\partial S$.
In what follows, a homotopy of essential arcs is always relative endpoints.

If $S$ is a surface with boundary,  we let $\mathcal{B}=\mathcal{B}(S)$ be the union of the set of homotopy classes
of essential arcs in $S$ with the set of homotopy classes of simple closed curves that are homotopic to boundary components. 
If $S$ is a surface without boundary, the set $\mathcal{B}$ is assumed to be empty.

For any surface $S$ with or without boundary, we consider the  function $J$ defined on $\mathcal{T}(S)\times \mathcal{T}(S)$ by 
$$
J(X,Y)=\log\sup_{\gamma\in\mathcal{C}\cup\mathcal{B}}\frac{l_{Y}(\gamma)}{l_{X}(\gamma)}
$$
for all $X,Y\in\mathcal{T}(S)$.
If the surface $S$ has no boundary, we recover Thurston's asymmetric metric $K$ defined above.

\begin{proposition}
The function $J: \mathcal{T}(S)\times \mathcal{T}(S)\to \mathbb{R}$ is an asymmetric metric on $\mathcal{T}(S)$.
\end{proposition}

\begin{proof}
The proof follows from \cite{LPST}, Propositions 2.10 and 2.13.
\end{proof}
 
It is shown in \cite{LPST}, Proposition 2.12, that when $S$ has nonempty boundary, the asymmetric metric
$J$ can be expressed as the logarithm of the supremum over the set $\mathcal{B}$ solely.   

In the same way as for the function $K$, we easily see that $J\leq L$.

\section{Surfaces of finite type}

We now construct Lipschitz-extremal homeomorphisms between some hyperbolic pairs of pants, using the homeomorphisms $h_k$ between symmetric hyperbolic right-angled hexagons that we constructed in Section \ref{s:extremal}. 
We shall then combine these homeomorphisms to get Lipschitz-extremal homeomorphisms of hyperbolic surfaces of arbitrary topological finite type.

We shall call a hyperbolic pair of pants {\it symmetric} if it is  obtained by gluing along three non-consecutive boundary components two isometric symmetric right-angled hexagons, and we shall always assume that these hexagons are glued along their long edges. 
Thus, the boundary components of our pairs of pants are ``short".

We let  $P$ be a symmetric pair of pants obtained by gluing two symmetric right-angled hexagons $H$, and for every $k\geq 0$, we let $P_k$ be a symmetric pair of pants obtained by gluing two right-angled hexagons $H_k$. 
Taking the double of the map $h_k:H\to H_k$ produces a map $p_k:P\to P_k$.

\begin{theorem} \label{th:geodesics}
The line $t\mapsto P_{e^{t}}$ ($t\in\mathbb{R}$) is a stretch line, and it is 
a geodesic for these two metrics $J$ and $L$ on $\mathcal{T}_{0,3}$. 
Furthermore, up to reparametrization, this line is also a geodesic for both metrics when it is traversed in the opposite direction. 
\end{theorem}

\begin{figure}[!hbp]
\psfrag{f}{$f_k$}
\psfrag{g}{$g_k$}
\psfrag{f}{$f_k$}
\psfrag{g}{$g_k$}
\psfrag{L}{$L$}
\psfrag{kL}{$kL$}
\psfrag{H}{$H$}
\psfrag{Hk}{$H_k$}
\centering
\includegraphics[width=0.6\linewidth]{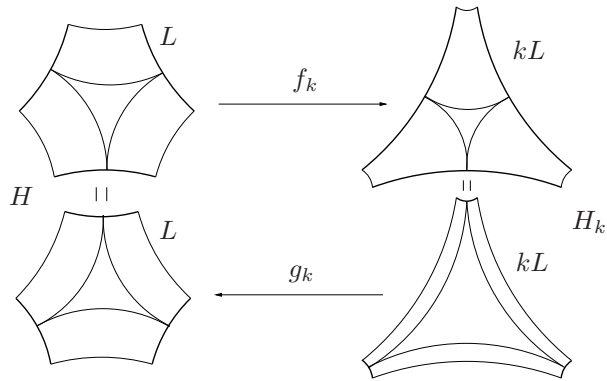}
\caption{\small{The actions of the maps $f_k$ and $g_k$ on symmetric hexagons.}}
\label{two}
\end{figure}

\begin{proof} 
For each $t \geq 0$, the action of the homeomorphism $P\to P_{e^{t}}$ on each boundary component of $P$ is linear (it multiplies arc length by $e^t$).
The fact that the line $t\mapsto P_{e^{t}}$ ($t\in\mathbb{R}$) coincides with a stretch line follows from the fact that for all $t\geq 0$, the surface $P_{e^{t}}$ is obtained from $P$ by multiplying the lengths of the boundary geodesics by the constant factor $e^t$, and this factor completely determines the resulting hyperbolic surface $P_{e^{t}}$. 
This also implies that we have $J(P,P_{e^{t}})=t$. 
On the other hand, since the map we construct is $e^t$-Lipschitz, we have $L(P,P_{e^{t}})\leq t$. 
This, together with the inequality $J\leq L$, gives $J(P,P_{e^{t}})=L(P,P_{e^{t}})$ for all $t\geq 0$. Thus, the line $t\mapsto P_{e^{t}}$  is a geodesic for $J$ and for $L$.

For the proof of the second statement, we first consider the case of hexagons. 
Let $H$ be a symmetric hexagon. 
Choose three non-consecutive edges as the long edges of $H$. 
For each $k\geq 1$, we have a map $h_k:H\to H_k$, as defined in Section \ref{s:extremal} above, whose Lipschitz constant is $k$ and which expands the long edges of $H$ by the factor $k$. 
By exchanging the roles of the long and short edges, we get a map $g_k:H_k\to H$ which expands the new long edges by a factor $d_k$, and contracts the new short edges  by the factor $k$.

From Formula (\ref{lengths}), we deduce that the dilatation factor $d_k$ of $g_k$ is given by
\begin{equation}\label{eq:dil}
d_k=\frac{l}{l_k}=
\displaystyle \frac
{\mathrm{argsinh \ }(\displaystyle\frac{1}{2\sinh L})}
 {\mathrm{argsinh \ }(\displaystyle\frac{1}{2\sinh kL})} 
\end{equation}

The homeomorphism $g_k$ has Lipschitz constant $d_k$ and it expands the long edges of the hyperbolic hexagon $H_k$ by the factor $d_k$ (see Figure \ref{two}), therefore we have
\[J(H_k,H)=\log d_k = L(H_k, H),\]
\[J(H,H_k)=\log k = L(H, H_k).\]
Doubling the hexagons, we obtain the same result for the symmetric pair of pants, showing that, up to parametrization, the line  $t\mapsto P_{e^{t}}$ is a geodesic in both directions for the metrics $J$ and $L$. 
\end{proof}
 
\begin{remark} 
By Theorem \ref{th:geodesics}, we have $J(x,y)=L(x,y)$ if the points $x$ and $y$ are situated on the stretch line that we construct.
We do not know whether the metrics $J$ and $L$ are equal on Teichm/"uller space. 
\end{remark}

A particular hyperbolic surface $S$ of arbitrary finite type $(g,b)$ can be obtained by gluing a collection of symmetric pairs of pants in such a way that the feet of the seams of adjacent pairs of pants coincide. In such a situation, we shall say that the gluing has been done {\it without torsion}.
For such a surface, we have the following:

\begin{theorem}\label{th:surfaces}
The line $t\mapsto S_{e^{t}}$ ($t\in\mathbb{R})$ is a stretch line in $\mathcal{T}_{g,b}$, and it is a geodesic for 
both asymmetric metrics $J$ and $L$ on $\mathcal{T}_{g,b}$. 
Up to reparametrization, this line  is also a geodesic for the same metrics when it is traversed in the opposite direction.
Along that line, the metrics $J$ and $L$ coincide. 
Furthermore, this stretch line has the following nice description in the Fenchel-Nielsen coordinates associated to the underlying pair of pants decomposition of $S$: 
at time $t$ from the origin, all the length parameters are multiplied by the constant factor $e^t$, and all the twist parameters are unchanged and remain equal to zero.
\end{theorem}

\begin{proof}
We start with a symmetric hyperbolic pair of pants $P$ equipped with a complete geodesic lamination, and we then consider the hyperbolic surface $S$, homeomorphic to $S_{0,4}$, obtained by gluing two copies of $P$ along one boundary component, in such a way that the following hold:
\begin{itemize}
\item The union of the complete geodesic laminations of both pairs of pants is a non chain-recurrent complete geodesic lamination of $S$.
\item The feet of the seams abutting on the component along which we glue coincide; that is, we glue without torsion. Here, the origin of Fenchel-Nielsen twist coordinates is measured as a signed distance between feets of seams (in the universal cover). We refer to  \cite[Theorem 4.6.23]{Thurston} for the convention on Fenchel-Nielsen coordinates.

\end{itemize} 
Let us denote by $\alpha$ the curve in $S$ that corresponds to the glued components.
There is an orientation-reversing order-two symmetry exchanging the copies of $P$ in $S$. 
The surface $S$ is equipped with a complete geodesic lamination $\mu$, and the order-two symmetry leaves the lamination $\mu$ invariant.

\begin{figure}[!hbp]
\psfrag{A}{$A$}
\psfrag{B}{$B$}
\centering
\includegraphics[width=0.4\linewidth]{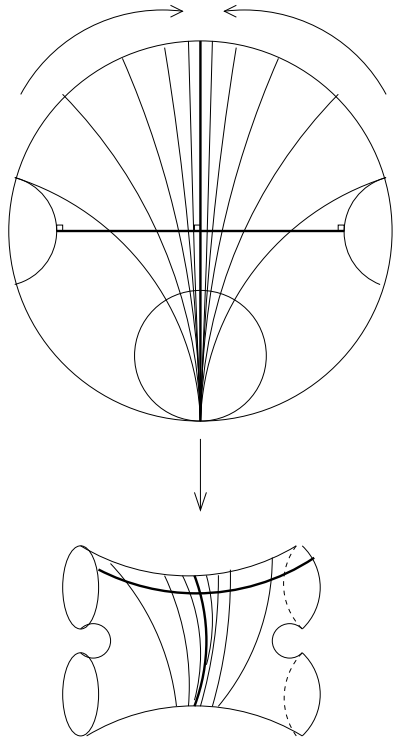}
\caption{\small{The action of a stretch map on the universal cover.}}
\label{haut}
\end{figure}

It is now useful to describe the situation in the universal covering $\widetilde{S}$ of $S$. 
The order-two symmetry lifts to the universal cover, and the preimage of $\mu$ in $\widetilde{S}$ is left invariant by this symmetry. 
The deformation of the hyperbolic plane by the stretch map can be seen in $\widetilde{S}$ as preserving a basepoint $O$ on a lift $\widetilde{\alpha}$ of $\alpha$ and the horocycle passsing through $O$ and centerd at the endpoint of $\widetilde{\alpha}$. 
The stretch deformation is then described in a neighborhood of $\widetilde{\alpha}$ by replacing the horocycle arcs that are contained in the spikes of each ideal triangle spiralling around $\widetilde{\alpha}$ by smaller arcs whose length has been raised to the power $e^t$. (Recall that the lengths of the horocycle pieces are all $<1$.) See Figure \ref{haut}) for a representation o this stretch deformation.
This shows that the stretch deformation commutes with the order-two symmetry. 
Hence, the feet of the seams coincide all along the deformation of $S$ by the stretch directed by $\mu$. 
In other words, stretching along $\mu$ does not induce Fenchel-Nielsen torsion.
The last statement of the theorem is thus established.
This also shows that the line $t\mapsto S_{e^{^t}}$ is a geodesic for both metrics $L$ and $J$, 
yielding the equality $L=J$ on that line.
We now proceed to show that our line traversed in opposite direction is a geodesic for both asymmetric metrics $J$ and $L$
and that these two metrics coincide along that line.
The homeomorphisms $g_{e^{t}}$ defined on each pair of pants given by the pants decomposition of $S$ piece together into a homeomorphism we also denote by $g_{e^{t}}$ from $S_{e^{t}}$ to $S$.
The reason why these local homeomorphisms piece together correctly is the absence of torsion along the components of the pants decomposition.
The Lipschitz constant of the homeomorphism $g_{e^{t}}$ thus obtained is $d_{e^{t}}$.
The seams of the pairs of pants coalesce into (smooth) geodesic simple closed curves and essential geodesic arcs that are stretched by the factor $d_{e^{t}}$ from $S_{e^{t}}$ to $S$.
This shows that the homeomorphism $g_{e^{t}}$ is Lipschitz-minimizing and that $L=J$ on the line.
The proof is complete.
\end{proof}
  
\begin{remark}  
The {\it dual metric} of an asymmetric metric $M$ on a set $X$ is the asymmetric metric defined by $\overline{M}(x,y)=M(y,x)$ for every $x$ and $y$ in $X$.
Equation (\ref{eq:dil}) shows that the asymmetric metric $J$ and its dual metric on $\mathcal{T}(S)$ are not quasi-isometric, even restricted to our geodesics $S_{e^t}$. 
Indeed, we have seen that for $t\geq0$, we have $J(S,S_{e^t})=t$ and $\overline{J}(S,S_{e^t})=\log d_{e^t}$.
But
 \begin{eqnarray*}
  d_{e^t}&\sim_{t\to\infty} &\mathrm{argsinh \ }\big(\displaystyle\frac{1}{2\sinh L}\big)\ \big/\  \mathrm{argsinh\ }\big(\displaystyle\frac{1}{2\sinh(e^{-t} L)}\big)\\
  &\sim_{t\to\infty} & \mathrm{argsinh \ }\big( \frac{1}{2\sinh L}\big)e^{e^tL},
 \end{eqnarray*}
that is, $\overline{J}(S,S_{e^t})\in O(e^t)$ as $t\to\infty$.
\end{remark}

Actually, we already noticed in \cite{PT1} and in \cite{Theret09} that Thurston's asymmetric metric for surfaces without boundary, of which $J$ is an analogue for surfaces with or without boundary, are not quasi-isometric to their dual metrics, in restriction to some special stretch lines. 
These observations naturally lead to the following:

\begin{question}
Characterize the geodesic lines for Thurston's asymmetric metric and for its analogue $J$ for surfaces with boundary, such that the restriction on that line of such a metric and its dual are quasi-isometric ?
\end{question}

We note in this respect that Choi and Rafi showed in \cite{CR} that in the thick part of Teichm\"uller space, Thurston's asymmetric metric and its dual metric are both quasi-isometric to Teichm\"uller's metric. 
On the other hand, there exist stretch lines that are completely contained in the thick part (take a pseudo-Anosov map whose stable and unstable laminations are complete, and consider the stretch line directed by one of these two laminations and passing by a point whose horocyclic foliation is the other lamination); 
therefore, there exist stretch lines for Thurston's asymmetric metric such that the restriction on that line of this metric and its dual are quasi-isometric.

We now recall that by a result of Thurston, given any two points $x$ and $y$ in Teichm\"uller space, there is a unique maximally stretched chain-recurrent geodesic lamination $\mu(x,y)$ from $x$ to $y$ which is maximal with respect to inclusion, and that if $x$ and $y$ lie in that order on a stretch line directed by a complete chain-recurrent geodesic lamination 
$\mu$, then $\mu(x,y)=\mu$. 
The next theorem identifies this geodesic lamination for two points $x$ and $y$ on the same stretch lines we construct, and it says in particular that this lamination is not complete.

\begin{theorem}\label{th:maxmax} For the stretch lines that we constructed above, 
the maximal maximally stretched lamination $\mu(S,S_{e^{t}})$ is the pair of pants decomposition that underlines the construction.
\end{theorem}

\begin{figure}[!hbp]
\centering
\psfrag{a}{\small \shortstack{horocycles\\perpendicular\\to the boundary}}
\psfrag{b}{\small \shortstack{horocyclic arc\\of length one}}
\psfrag{c}{\small \shortstack{non-foliated\\region}}
\includegraphics[width=0.25\linewidth]{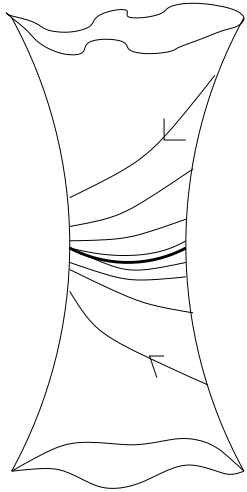}
\caption{\small {A non chain-recurrent geodesic lamination. The spirals wrap around the closed curve in opposite directions.}}
\label{spiral}
\end{figure}

\begin{proof}
Let $t>0$. 
The maximal maximally stretched chain-recurrent geodesic lamination $\mu(S,S_{e^{t}})$ from $S$ to $S_{e^{t}}$ contains the underlying pair of pants decomposition, since each curve in this decomposition is maximally stretched. 
Assume for contradiction that  $\mu(S,S_{e^{t}})$ contains a larger lamination. 
It then contains a bi-infinite geodesic that spirals around some closed geodesic $C$ in that decomposition. 
Since $\mu(S,S_{e^{t}})$ is chain-recurrent, it contains another geodesic that spirals along the opposite side of $C$ in the same direction (compare Figure \ref{spiral}). 
By a result in \cite{PT2}, if we perform a Thurston stretch along a completion of $\mu(S,S_{e^{t}})$, then we necessarily introduce a Fenchel-Nielsen torsion about the closed geodesic $C$. Now
Thurston proved in \cite{Thurston-stretch} that we can join $S$ to $S_{e^t}$ by a concatenation of Thurston stretches which are directed by complete geodesic laminations, all of them containing $\mu(S,S_{e^{t}})$. 
The torsions introduced about the geodesic $C$ are all in the same direction. 
Thus, there necessarily is a nonzero torsion. 
This contradicts Theorem \ref{th:surfaces}. 
Thus, $\mu(S,S_{e^{t}})$ does not contain any geodesic  lamination larger than the geodesics of the pair of pants decomposition.
Thus, the maximal maximally stretched lamination $\mu(S,S_{e^{t}})$ is the pair of pants decomposition.
\end{proof}

It also follows from the reasoning in the proof of Theorem \ref{th:maxmax} that the set of maximal maximally stretched laminations from $S$ to $S_{e^{t}}$ is the set of all completions of the pants decomposition that are nowhere chain-recurrent, which means that the geodesics spiralling around each component of the pants decomposition wrap in opposite directions, as illustrated in Figure \ref{spiral}.

\begin{figure}[!hbp]
\centering
\psfrag{a}{\small \shortstack{horocycles\\perpendicular\\to the boundary}}
\psfrag{b}{\small \shortstack{horocyclic arc\\of length one}}
\psfrag{c}{\small \shortstack{non-foliated\\region}}
\includegraphics[width=0.45\linewidth]{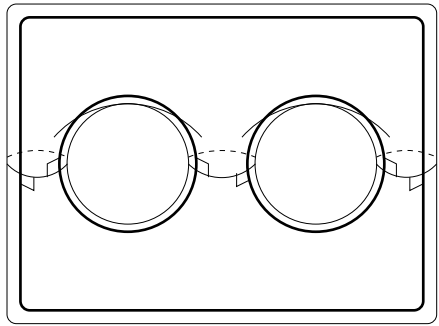}
\caption{\small {In bold lines is represented a  pants decompostion of the closed surface of genus $2$. 
The union of the seams is a multi-curve and a pants decomposition as well for the genus 2 surface.}}
\label{multicurve}
\end{figure}

\begin{remark} 
Given two points $x,y$ in Teichm\"uller space and knowing the maximal maximally stretched lamination $\mu(x,y)$ 
from $x$ to $y$, it is in general quite difficult to find the lamination $\mu(y,x)$.
For all $t>0$,  the maximal maximally stretched "lamination" from $S_{e^{t}}$ to $S$ is the union of the seams. 
As already mentioned in the proof of Theorem \ref{th:surfaces}, by our choice of the twist parameters (in which the feet of the seams coincide), in the case of closed surfaces, the union of the seams is a union of disjoint closed geodesics (a multi-curve), see Figure \ref{multicurve}. 
This multi-curve is maximally stretched by the stretch that we defined from $S_{e^{t}}$ to $S$ and therefore it is contained in the lamination $\mu(S_{e^{t}},S)$. 
In the case of a closed surface of genus $2$, the preceding argument shows that $\mu(S_{e^{t}},S)$ \emph{is} a union of seams, since this union is a pants decomposition.
\end{remark}

\end{document}